\documentclass[aos,preprint]{imsart}

\RequirePackage{amsthm,amsmath,amsfonts,amssymb,bm,graphicx,booktabs,
  multirow,threeparttable,mathtools,placeins,array}
\RequirePackage[round,authoryear]{natbib}
\RequirePackage{xcolor}

\RequirePackage[hidelinks]{hyperref}

\startlocaldefs
\newtheorem{theorem}{Theorem}
\newtheorem{proposition}{Proposition}
\newtheorem{lemma}{Lemma}

\theoremstyle{definition}

\newcommand{\R}{\mathbb{R}}
\newcommand{\E}{\mathbb{E}}
\newcommand{\Var}{\mathrm{Var}}
\newcommand{\Cov}{\mathrm{Cov}}

\newcommand{\Prb}{\mathbb{P}}
\newcommand{\N}{\mathcal{N}}
\newcommand{\tr}{\operatorname{tr}}
\newcommand{\1}{\mathbf{1}}
 \endlocaldefs

\begin{document}

\begin{frontmatter}
\title{Dimension Changes the Scale-Signal Channel in
Nearest-Neighbor Two-Sample Testing}
\runtitle{Dimension Changes the Scale-Signal Channel}
\begin{aug}
\author[A]{\fnms{Hao}~\snm{Chen}\ead[label=e1]{hxchen@ucdavis.edu}}
\address[A]{Department of Statistics, University of California, Davis
\printead[presep={,\ }]{e1}}
\end{aug}

\begin{abstract}
Nearest-neighbor two-sample tests are often analyzed with dimension fixed. We study the power of the generalized edge-count test on a directed 1-nearest-neighbor graph. This statistic combines a weighted sum and a difference of the two within-sample edge counts. For balanced isotropic Gaussian samples with $d/\log N\to\gamma>0$, we derive the connection probabilities under local location and scale alternatives. Under location alternatives, the two counts have equal mean shifts, which cancel in their difference. Under scale alternatives, the two counts instead have mean shifts of opposite signs. We also derive the growth rate of the null expectation $\E_0C_3$, where $C_3$ counts pairs of directed edges sharing exactly one vertex. On an explicit range of $\gamma$, we establish local-power limits strictly above the test level for scale alternatives in the difference component at the scale $|\varepsilon_N|\asymp\sqrt{\E_0C_3}/(Nd)$. Here $\E_0C_3$ gives the variance order of the difference count under these local alternatives. At this scale, the standardized drift of the weighted sum vanishes. When $d=1$, the weighted component instead detects local Gaussian scale alternatives while the difference component has no limiting drift. Thus dimension changes both the magnitude of a detectable scale difference and the component of the test that detects it.
\end{abstract}

\begin{keyword}[class=MSC]
\kwd{Primary 62G10; secondary 62G20, 62G09}
\end{keyword}

\begin{keyword}
\kwd{two-sample test}
\kwd{nearest-neighbor graph}
\kwd{edge-count test}
\kwd{high-dimensional inference}
\kwd{local power}
\end{keyword}
\end{frontmatter}
\hypersetup{bookmarksdepth=2,
pdfkeywords={two-sample test, nearest-neighbor graph, edge-count test, high-dimensional inference, local power}}

\section{Introduction}
\label{sec:intro}

Let $X_1,\ldots,X_m\overset{\mathrm{iid}}{\sim}F_X$ and $Y_1,\ldots,Y_n\overset{\mathrm{iid}}{\sim}F_Y$ be independent samples in $\R^d$, and write $N=m+n$.  We test $H_0:F_X=F_Y$ against $H_a:F_X\ne F_Y$.  Graph-based two-sample tests pool the observations and construct a similarity graph without using the sample labels. Because this construction requires only a dissimilarity, it also applies to high-dimensional and non-Euclidean data. Classical edge-count tests compare the number of between-sample edges with its permutation distribution \citep{friedman1979,schilling1986,henze1988,rosenbaum2005}. \citet{chen2017} introduced the generalized edge-count test, which jointly uses the two within-sample edge counts. This statistic can be decomposed into weighted and difference components \citep{chen2018,chu2019}.  The standardized weighted sum, $Z_w$, measures the overall excess of within-sample connections, while the standardized difference, $Z_d$, contrasts the within-sample connections of the two samples.

Nearest-neighbor two-sample tests are usually studied with dimension fixed. For fixed $d$ and a fixed number $k$ of nearest neighbors, these tests have zero Pitman efficiency against $N^{-1/2}$ parametric alternatives \citep{bhattacharya2019}. \citet{bhattacharya2020} characterizes the detection threshold for different fixed dimensions when $k$ is fixed, while \citet{kanekar2026} studies how it changes when $k$ grows with sample size.

In this paper, we keep one outgoing edge per observation and let dimension grow with sample size. For balanced isotropic Gaussian samples with $d/\log N\to\gamma>0$, we show that the two within-sample connection probabilities respond in opposite directions at first order under scale alternatives. On an explicit range of $\gamma$, a local relative scale change $\varepsilon_N$ produces a mean shift of order $Nd|\varepsilon_N|$ in the difference count, whose variance is of order $\E_0C_3$. Here $C_3$ counts pairs of directed edges sharing exactly one vertex. Comparing this mean shift with the standard deviation gives the local scale
\[
 |\varepsilon_N|\asymp\frac{\sqrt{\E_0C_3}}{Nd}.
\]
We establish nondegenerate local-power limits on this scale. If $\varepsilon_NNd/\sqrt{\E_0C_3}$ converges to a nonzero constant, the difference component $Z_d$ has a nonzero limiting standardized drift, whereas that of the weighted component $Z_w$ vanishes. We also show that $\E_0C_3=N^{1+{\beta_\gamma}+o(1)}$, where ${\beta_\gamma}>0$ is given by the variational formula in Section~\ref{sec:permclt}. Thus the local scale above is $N^{-1/2+{\beta_\gamma}/2+o(1)}/d$. For example, at $\gamma=0.25$, numerical evaluation gives ${\beta_\gamma}\approx1.4\times10^{-8}$, yielding a local scale of approximately $N^{-0.499999993+o(1)}/d$. For comparison, we prove that in one dimension, scale alternatives $\varepsilon_N=aN^{-1/4}$ with fixed $a\ne0$ produce nonzero limiting drift in $Z_w$ and zero limiting drift in $Z_d$.

For balanced Gaussian location alternatives, the expected difference between the two within-sample edge counts remains zero, and the standardized weighted sum $Z_w$ has nondegenerate local power at $\|\mu_N\|\asymp N^{-1/4}$, both in one dimension and in the logarithmic dimension range considered here.

The power analysis has three parts.  First, exact conditional integrals reduce the connection probabilities to the minima of within-sample and between-sample distances.  A conditional strong-large-deviation expansion describes the noncentral chi-square lower tails in logarithmic dimension, and a Gaussian conditional calculation gives the expected squared norm of the selected nearest neighbor that controls the scale derivative.  Second, exact permutation formulas express the two component variances in terms of mutual-neighbor and shared-vertex counts.  In logarithmic dimension, the fraction of mutual-neighbor edges vanishes and ${\E_0}C_3$ grows superlinearly.  Third, we compare the sampling and permutation covariances and obtain normal approximations under both distributions. Combining these results with the changes in the expected edge counts under local alternatives gives the local-power limits in logarithmic dimension.  The one-dimensional power result combines localization and density-overlap calculations with stabilization and convergence of the permutation critical value.

\begin{figure}[!tb]
\centering
\includegraphics[width=0.88\textwidth]{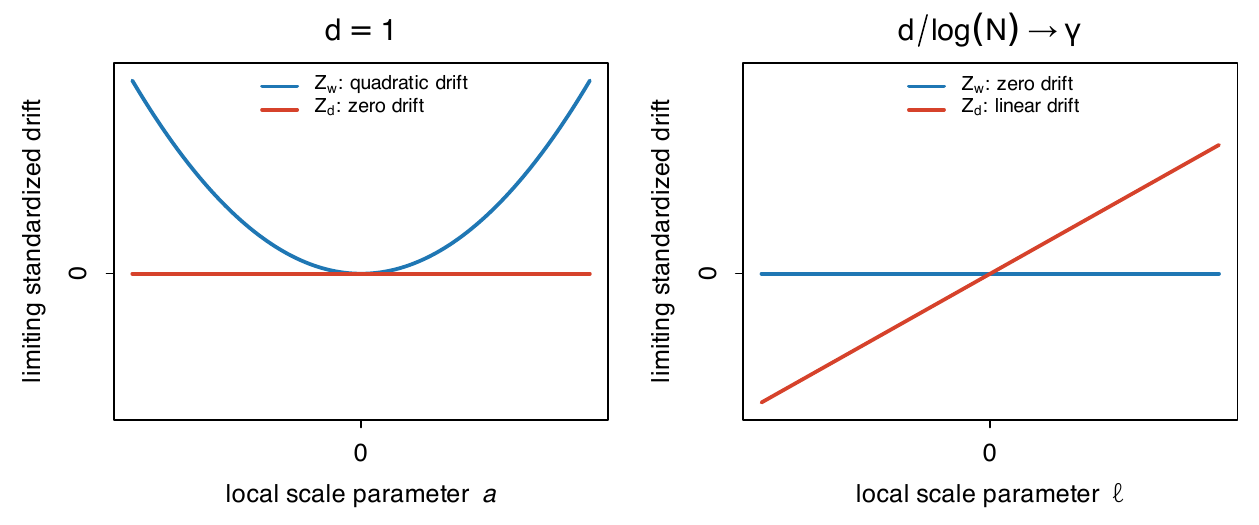}
\caption{Limiting standardized drifts under balanced isotropic Gaussian scale alternatives. Left: $d=1$, with $a=N^{1/4}\varepsilon_N$. Right: $d/\log N\to\gamma$ in the dimension range of Theorem~\ref{thm:logscalepower}, with $\ell=\varepsilon_NNd/\sqrt{\E_0C_3}$. Schematic curves; model-specific positive constants are suppressed.}
\label{fig:channelendpoints}
\end{figure}

Figure~\ref{fig:channelendpoints} summarizes these two scale-power results. In one dimension, the limiting drift of $Z_w$ is quadratic in $a=N^{1/4}\varepsilon_N$, and that of $Z_d$ is zero (Theorem~\ref{thm:fixedendtoend}). In logarithmic dimension, the limiting drift of $Z_d$ is linear in $\ell=\varepsilon_NNd/\sqrt{\E_0C_3}$, and that of $Z_w$ is zero (Theorem~\ref{thm:logscalepower}).
\FloatBarrier

\paragraph*{Related work} \citet{henze1999} established universal consistency of the \linebreak Friedman--Rafsky MST test in fixed dimension.  A general local-power theory for graph-based tests, including the Chen--Friedman statistic, was developed by \citet{bhattacharya2019}; it also shows how graph density affects Pitman efficiency. \citet{bhattacharya2020} established normal limits for stabilizing graph statistics in a Poissonized setting and derived detection thresholds for the directed $k$-NN count of edges from the first sample to the second. The formal parametric results assume compact support. The analysis also covers the total between-sample edge count, including both directions. For this count, the first-order coefficient vanishes when the two sample sizes are equal. \citet{kanekar2026} studies the distinct fixed-dimensional problem in which the number of nearest neighbors grows.

Newman, Rinott and Tversky \citeyearpar{newman1983} proved that the nearest-neighbor in-degree variance diverges when dimension tends to infinity before sample size. \citet{schilling1986mutual} derived mutual- and shared-neighbor probability limits in fixed dimension and then studied their behavior as dimension tends to infinity. Nearest-neighbor hubs in high dimension were documented by \citet{radovanovic2010}. High-dimension, low-sample-size graph and distance tests are studied by \citet{biswas2014}, \citet{mondal2015}, and \citet{li2018}; those asymptotics hold sample size fixed or use central limits for interpoint distances.  \citet{sarkar2020} modify the dissimilarity used to construct graph-based tests and establish consistency as $d\to\infty$ with $m$ and $n$ fixed.  Power analyses with both dimension and sample size tending to infinity have been developed for kernel two-sample statistics \citep{reddi2015,yan2023}.

\section{Test statistics and connection probabilities}
\subsection{The test statistic and its components}
\label{sec:setup}

We first define the statistic. Pool the observations as $Z_1,\ldots,Z_N$, and let $g_i=1$ for sample $X$ and $g_i=0$ for sample $Y$.  Construct the directed 1-NN graph $G$ by drawing an edge from each observation to its nearest neighbor. Define
\[
R_1=\sum_{(i,j)\in G}\1\{g_i=g_j=1\},\qquad
R_2=\sum_{(i,j)\in G}\1\{g_i=g_j=0\}.
\]
Let $R=(R_1,R_2)^\top$.  The subscript $P$ denotes the permutation distribution conditional on the pooled graph, with the sample sizes $m$ and $n$ fixed. Write $\mu_P$ and $\Sigma_P$ for the exact permutation mean and covariance of $R$.  The directed generalized edge-count statistic is\footnote{For $m,n\ge2$ and unique nearest neighbors, $\Sigma_P$ is singular if and only if every vertex has in-degree $d_i^{\mathrm{in}}=1$. In that case, the inverse is taken on the subspace with positive permutation variance.}
\begin{equation}
S=(R-\mu_P)^\top{\Sigma_P^{-1}}(R-\mu_P).
\label{eq:decomp}
\end{equation}
Define the original, weighted, and difference counts by
\begin{equation}
 R_o={|G|-R_1-R_2},
 \qquad
 R_w=\frac{n-1}{N-2}R_1+\frac{m-1}{N-2}R_2,
 \qquad
 R_d=R_1-R_2.
 \label{eq:originalweighteddefinitions}
\end{equation}
Here ${|G|}$ is the number of directed edges, and $R_o$ counts between-sample edges. Write $\sigma_j^2=\Var_P(R_j)$ and put
\[
 {Z_j=\frac{R_j-\E_PR_j}{\sigma_j},\qquad j\in\{o,w,d\},\quad \sigma_j^2>0.}
\]
Under the permutation distribution, $\Cov_P(R_w,R_d)=0$, and the statistic decomposes as $S=Z_w^2+Z_d^2$ \citep{liu2022}.

Let $c_{\alpha,N}(G)$ be the exact permutation critical value.  The test rejects when $S>c_{\alpha,N}(G)$.  Under the null hypothesis, the exact permutation test has rejection probability at most $\alpha$ for every sample size. The discreteness of the permutation distribution can make the test conservative.

The original edge count $R_o$ can also be expressed in terms of the weighted and difference counts $R_w$ and $R_d$.  For every directed graph without self-loops and every $2\le m,n\le N-2$,
\begin{equation}
 R_o-\E_PR_o
 ={-}2(R_w-\E_PR_w)\mathbin{-}\frac{m-n}{N-2}(R_d-\E_PR_d).
 \label{eq:originalmixture}
\end{equation}
Then
\begin{equation}
 Z_o={-}\frac{Z_w+r_NZ_d}{\sqrt{1+r_N^2}},
 \qquad{\text{where}}\qquad
 r_N=\frac{m-n}{N-2}\frac{\sigma_d}{2\sigma_w}.
 \label{eq:originalstandardizedmixture}
\end{equation}
In particular, at balance $R_o=|G|-2R_w$, so $Z_o={-}Z_w$. Thus, when $m=n$, the classical test based on small values of $Z_o$ is equivalent to the test based on large values of $Z_w$.

{For balanced location alternatives, central symmetry forces the two expected within-sample counts to agree at every finite sample size.} Let $m=n$, let $F$ be centrally symmetric on $\mathbb R^d$ {with its center taken to be zero without loss of generality}, and suppose
\[
 X_i\stackrel{\rm iid}{\sim}F,\qquad
 Y_j\stackrel{\rm iid}{\sim}\mu+F.
\]
Suppose the graph is constructed without using sample labels, and its connections are unchanged by translations, rotations, reflections, or reordering of the observations.\footnote{Any ties are resolved in a way that preserves these properties.} Then
\begin{equation}
 (R_1,R_2)^{\top}\overset d=(R_2,R_1)^{\top}.
 \label{eq:exactlocationsymmetry}
\end{equation}
Consequently, under the location alternative,
\begin{equation}
 \E_{a} R_d=0,
 \qquad
 \Cov_{a}(R_w,R_d)=0{,}
 \label{eq:exactlocationchannels}
\end{equation}
where the subscript $a$ denotes the sampling distribution under the alternative. This applies, in particular, to the directed 1-NN graph at every finite $(N,d)$.  By central symmetry, $\mu-X_i$ has distribution $\mu+F$, and $\mu-Y_j$ has distribution $F$. Reflecting all observations by $z\mapsto\mu-z$ preserves their pairwise distances. Exchanging the two samples then preserves the joint distribution of the data and interchanges $R_1$ and $R_2$, giving \eqref{eq:exactlocationsymmetry}. Since $m=n$, $R_w=(R_1+R_2)/2$ and $R_d=R_1-R_2$. Hence
\[
 {\E_a R_d=0,\qquad
 \Cov_a(R_w,R_d)=\frac12\{\Var_a(R_1)-\Var_a(R_2)\}=0.}
\]

\subsection{Connection probabilities}
\label{sec:exact}
\label{sec:gaussian}

Let $p_1$ be the probability that a randomly selected observation from sample $X$ has its nearest neighbor in the same sample, and define $p_2$ analogously for sample $Y$. The expected within-sample edge counts under the alternative satisfy
\begin{equation}
\E_aR_1=mp_1,\qquad \E_aR_2=np_2.
\label{eq:means}
\end{equation}

A nearest neighbor belongs to the same sample when the smallest within-sample distance is less than the smallest between-sample distance.  We first give an exact formula for this probability, then identify the relevant lower tails and evaluate them for Gaussian samples.

Condition on an observation $x$ from sample $X$.  Let $F_{XX,x}$ be the conditional distribution function of $\|x-X\|$ for an independent $X\sim F_X$, and define $F_{XY,x}$ in the same way using $Y\sim F_Y$.  Assume that the distance distributions are continuous, and write $n_w=m-1$ and $n_b=n$. Let $p_1(x)$ be the probability that an observation from sample $X$ at $x$ has its nearest neighbor in sample $X$. The finite-sample connection formula of \citet[p.~801, proof of Theorem~3.4]{schilling1986}, written conditionally on $x$, gives
\begin{align}
p_1(x)
&=n_w\int [1-F_{XX,x}(t)]^{n_w-1}[1-F_{XY,x}(t)]^{n_b}
\,dF_{XX,x}(t),\label{eq:p1exact}\\
p_1&=\E_{X_0\sim F_X}\{p_1(X_0)\}.\nonumber
\end{align}
The formula for $p_2$ follows by interchanging $X$ and $Y$ and using $n_w=n-1$ and $n_b=m$.

Fix an observation $x$ from sample $X$ and one of the other $m-1$ observations from that sample. If their distance is $t$, this candidate is the nearest neighbor when all remaining distances from $x$ exceed $t$. Conditional on $x$, the remaining observations are independent, so this event has probability $[1-F_{XX,x}(t)]^{n_w-1}[1-F_{XY,x}(t)]^{n_b}$. Integrating over $t$ and summing over the $n_w=m-1$ candidates gives $p_1(x)$. Averaging over $x$ gives $p_1$.

Given $x$, the probability that all $n_w$ within-sample distances exceed $t$ is $[1-F_{XX,x}(t)]^{n_w}$. As $n_w$ grows, this probability stays away from both zero and one when $F_{XX,x}(t)$ is of order $1/n_w$. Given $x$, the probability that all $n_b$ between-sample distances exceed $t$ is $[1-F_{XY,x}(t)]^{n_b}$. As $n_b$ grows, this probability stays away from both zero and one when $F_{XY,x}(t)$ is of order $1/n_b$. These are the probability levels relevant to comparing the two minima.

We now obtain these conditional distance distributions for Gaussian samples. Let $F_X=\N(0,I_d)$ and $F_Y=\N(\mu,\Sigma)$, where $\Sigma\succ0$. Let ${Z_1}$ be an observation from either sample, and let $Z\sim\N(a,C)$ be another observation, independent of ${Z_1}$. Here $(a,C)=(0,I_d)$ for a candidate from sample $X$, and $(a,C)=(\mu,\Sigma)$ for one from sample $Y$.

Conditional on ${Z_1}={z_1}$, write $Q=\|{z_1}-Z\|^2$ for their squared distance. Write $C=U\operatorname{diag}(\lambda_1,\ldots,\lambda_d)U^\top$ and put $\theta=U^\top({z_1}-a)$. The vector $U^\top(Z-a)$ is Gaussian with mean zero and covariance $\operatorname{diag}(\lambda_1,\ldots,\lambda_d)$. Hence its coordinates are independent, and
\[
 {W_j=\frac{[U^\top(Z-a)]_j}{\sqrt{\lambda_j}},\qquad j=1,\ldots,d,\qquad W_j\overset{\rm iid}{\sim}\N(0,1).}
\]
Then
\begin{align*}
Q&=\|{z_1}-Z\|^2\mathrel{=}\|U^\top({z_1}-Z)\|^2
={\textstyle\sum_{j=1}^d}\left(\theta_j-\sqrt{\lambda_j}W_j\right)^2,\\
\E Q&\mathrel{=}{\textstyle\sum_{j=1}^d}(\theta_j^2+\lambda_j)
=\|{z_1}-a\|^2+\tr C,\\
\Var Q&\mathrel{=}{\textstyle\sum_{j=1}^d}\left(2\lambda_j^2+4\lambda_j\theta_j^2\right)
=4({z_1}-a)^\top C({z_1}-a)+2\tr C^2.
\end{align*}

In the isotropic case $C=\sigma^2I_d$, this gives the exact conditional law
\[
 {\frac{Q}{\sigma^2}\ \sim\ \chi_d^2\!\left(\frac{\|{z_1}-a\|^2}{\sigma^2}\right).}
\]
For an observation ${z_1}=x$ from sample $X$, the two choices of $(a,C)$ give $F_{XX,x}(t)$ and $F_{XY,x}(t)$ through $\Prb(Q\le t^2\mid {Z_1}=x)$, for $t\ge0$. The same calculation applies to an observation from sample $Y$.

When $d/\log N\to\gamma>0$, we use large-deviation expansions of these conditional laws to calculate how the connection probabilities respond to location and scale changes, with the resulting mean shifts and power limits given in Section~\ref{sec:dichotomy}. The conditional-law Monte Carlo calculations in Section~\ref{sec:sims} use these exact conditional laws.

\section{Graph geometry and permutation distribution}
\label{sec:permclt}

We next determine how graph geometry affects the variances and covariance of $R_1$ and $R_2$. We then use a normal approximation under the permutation distribution to obtain the permutation critical value. Section~\ref{sec:dichotomy} combines these results with mean shifts and a normal approximation under the alternative to derive the local-power limits.

\subsection{Graph geometry and component variances}
\label{sec:permutationmoments}

Conditioning on the pooled graph, the permutation variances and covariance of $R_1$ and $R_2$ can be expressed in terms of two edge-pair counts. Let $C_2$ and $C_3$ count the ordered pairs of distinct directed edges that are opposite or share exactly one vertex, respectively. For directed $k$-NN graphs, the exact permutation means, variances, and covariance are given by \citet[Theorem~1]{liu2022}. The means are
\[
\E_PR_1={kN}{\frac{m(m-1)}{N(N-1)}},\qquad \E_PR_2={kN}{\frac{n(n-1)}{N(N-1)}}.
\]
The Supplement gives the full variance and covariance formulas in our notation.

Let $m=n=N/2$ and use the directed 1-NN graph. These exact moments give
\begin{equation}
\begin{aligned}
 {\Var_P(R_w)}
 &{=\frac{N}{16(N-1)(N-3)}
 \left\{N(N-2)+(N-4)C_2-C_3+\frac{2N}{N-1}\right\},}\\
 {\Var_P(R_d)}&{=\frac{N}{4(N-1)}(C_3+2C_2-2N).}
\end{aligned}
 \label{eq:balancedcomponentvariances}
\end{equation}
These formulas reduce the variance analysis to the behavior of $C_2$ and $C_3$.

We first give an identity relating $C_2$, $C_3$, and the vertex degrees. For vertex $i$, let $d_i^{\rm out}$ and $d_i^{\rm in}$ be the numbers of edges leaving and entering $i$, respectively, and let ${d_i}=d_i^{\rm out}+d_i^{\rm in}$ be the total directed degree of vertex $i$.  For the directed 1-NN graph, ${d_i}=1+{d_i^{\rm in}}$. At vertex $i$, there are ${d_i}({d_i}-1)$ ordered pairs of distinct directed edges sharing that vertex. Summing over vertices gives
\begin{equation}
 C_3=\sum_{i=1}^N{d_i}({d_i}-1)-2C_2.
\label{eq:c3degree}
\end{equation}
Each ordered pair of opposite edges is counted at both endpoints and must be excluded from $C_3$, giving the subtraction of $2C_2$.

Write $\E_0$ for expectation over the data under the null hypothesis. We next calculate the probability that two observations have the same nearest neighbor, which determines the expected number of pairs of incoming edges.  Let ${Z_1},{Z_2},{Z_3},Z$ be iid from the pooled null distribution, and put
\[
 {q_{1;23}}=\Prb\bigl(\|Z-{Z_2}\|^2\mathrel{>} {\|{Z_1}-{Z_2}\|^2}{,}\
 \|Z-{Z_3}\|^2\mathrel{>} {\|{Z_1}-{Z_3}\|^2}\mid {Z_1},{Z_2},{Z_3}\bigr).
\]
Given ${Z_1},{Z_2},{Z_3}$, this is the probability that another independent observation $Z$ is farther from both ${Z_2}$ and ${Z_3}$ than ${Z_1}$ is. For ${Z_1}$ to be the nearest neighbor of both ${Z_2}$ and ${Z_3}$, both observations must be closer to ${Z_1}$ than to each other, that is,
\[
 {\|{Z_2}-{Z_3}\|^2>\max\{\|{Z_1}-{Z_2}\|^2,\|{Z_1}-{Z_3}\|^2\}.}
\]
Conditional on ${Z_1},{Z_2},{Z_3}$, each of the remaining $N-3$ observations lies outside both balls centered at ${Z_2}$ and ${Z_3}$, with radii $\|{Z_1}-{Z_2}\|$ and $\|{Z_1}-{Z_3}\|$, with probability ${q_{1;23}}$. Conditional independence gives the factor ${q_{1;23}^{N-3}}$. Thus, if the distance condition above holds, the conditional probability that both ${Z_2}$ and ${Z_3}$ have ${Z_1}$ as their nearest neighbor is ${q_{1;23}^{N-3}}$; otherwise, it is zero.

For each vertex $i$, $d_i^{\rm in}(d_i^{\rm in}-1)$ counts ordered pairs of distinct observations whose nearest neighbor is $i$. There are $N(N-1)(N-2)$ ordered triples of distinct observations. By exchangeability and the conditional probability above,
\[
{\frac1N\E_0\sum_{i=1}^N d_i^{\rm in}(d_i^{\rm in}-1)
 =(N-1)(N-2)\E_0[{q_{1;23}^{N-3}}\1_{\{\|{Z_2}-{Z_3}\|^2>\max(\|{Z_1}-{Z_2}\|^2,\|{Z_1}-{Z_3}\|^2)\}}].}
\]
Since ${d_i}=1+d_i^{\rm in}$ and $\sum_i d_i^{\rm in}=N$,
\[
 {\sum_{i=1}^N{d_i}({d_i}-1)
 =\sum_{i=1}^N\left\{d_i^{\rm in}(d_i^{\rm in}-1)+2d_i^{\rm in}\right\}
 =\sum_{i=1}^N d_i^{\rm in}(d_i^{\rm in}-1)+2N.}
\]
Substituting into \eqref{eq:c3degree}, taking expectations, and dividing by $N$ gives
\begin{equation}
\begin{aligned}
 \frac1N\E_0C_3
 ={}&2+(N-1)(N-2)\E_0\left[
 {{q_{1;23}^{N-3}}}
 {\1_{\{\|{{Z_2}}-{{Z_3}}\|^2>\max(\|{{Z_1}}-{{Z_2}}\|^2,\|{{Z_1}}-{{Z_3}}\|^2)\}}}
 \right]\\
 &-\frac2N\E_0C_2.
\end{aligned}
\label{eq:c3threesourcemain}
\end{equation}

We now study these graph counts under the isotropic Gaussian null when $d/\log N\to\gamma>0$, beginning with the geometry of an observation and its nearest neighbor. Theorem~\ref{thm:logendpointgeometry} determines the normalized squared-norm limit of the selected nearest neighbor, which we use below to show that $C_2/N\to0$. It also gives the limit of the nearest squared distance divided by $d$, using the conditional Gaussian squared-distance tail analyzed by \citet[Section~IV-D]{zhou2019} for minimum-distance selection from an i.i.d. Gaussian codebook.  Let $s_\gamma\in(0,1)$ be the unique solution of
\begin{equation}
 \frac12\{s_\gamma^2-s_\gamma-\log s_\gamma\}=\gamma^{-1},
 \qquad y_\gamma=s_\gamma+s_\gamma^2.
\label{eq:gammasaddle}
\end{equation}

\begin{theorem}
\label{thm:logendpointgeometry}
Let ${N\to\infty}$, and let ${Z_1},{Z_2},\ldots,{Z_N}$ be independent $\N(0,I_d)$ vectors.  Set
\[
 {i_*}=\arg\min_{2\le i\le N}\|Z_i-{Z_1}\|.
\]
If $d\to\infty$ and $d/\log N\to\gamma\in(0,\infty)$, then
\begin{equation}
 \frac1d\min_{2\le i\le N}\|Z_i-{Z_1}\|^2
 \xrightarrow{L^1}y_\gamma,
 \qquad
 \frac1d\|Z_{i_*}\|^2\mathrel{\xrightarrow{L^1}}
 1-s_\gamma+s_\gamma^2.
 \label{eq:logendpointnorm}
\end{equation}
Consequently,
\begin{equation}
 d-\E_0\|Z_{i_*}\|^2
 \sim d\,s_\gamma(1-s_\gamma).
 \label{eq:logendpointdeficit}
\end{equation}
The convergence is uniform when $\gamma$ ranges over a compact subset of $(0,\infty)$.
\end{theorem}

Although all observations have the same marginal distribution, nearest-neighbor selection favors observations closer to the distribution mean: the selected neighbor's limiting normalized squared norm is below the typical value $1$. In contrast, in fixed dimension, the nearest neighbor approaches $Z_1$ as $N\to\infty$, and $\|Z_{i_*}\|^2-\|Z_1\|^2\to0$ in $L^1$.

{The proof of Theorem~\ref{thm:logendpointgeometry} first conditions on $Z_1$ and uses a noncentral chi-square lower-tail expansion to {show that $\|Z_{i_*}-Z_1\|^2/d\to y_\gamma$ in probability}.} To obtain the neighbor's norm, condition further on $Z_1=x$ and $\|Z_{i_*}-x\|^2=w$. Selecting the smallest distance does not change the conditional distribution of the direction at that distance. If $u$ is the cosine of the angle between $x$ and $Z_{i_*}-x$, its conditional density is proportional to
\[
 {\exp(-\|x\|\sqrt w\,u)(1-u^2)^{(d-3)/2},\qquad -1<u<1.}
\]
Together with $\|Z_1\|^2/d\to1$, this distance limit implies that the logarithm of this density divided by $d$, up to an additive constant, approaches $-\sqrt{y_\gamma}\,u+\tfrac12\log(1-u^2)$. Its unique maximum is at $u=-s_\gamma/\sqrt{y_\gamma}$. A Laplace approximation therefore gives $Z_1^\top(Z_{i_*}-Z_1)/d\to-s_\gamma$ in probability, and hence
\[
 {\begin{aligned}
 \frac{\|Z_{i_*}\|^2}{d}
 &=\frac{\|Z_1\|^2}{d}+\frac{\|Z_{i_*}-Z_1\|^2}{d}
   +\frac{2Z_1^\top(Z_{i_*}-Z_1)}{d}\\
 &\xrightarrow{p}1+y_\gamma-2s_\gamma=1-s_\gamma+s_\gamma^2.
 \end{aligned}}
\]
Moment bounds give uniform integrability and hence convergence in $L^1$. The detailed proof is given in the Supplement.

To determine the growth of ${\E_0}C_3$, we use the shared-neighbor probability in \eqref{eq:c3threesourcemain}. Under the isotropic Gaussian null, the squared distance from a fixed observation to an independent candidate has a noncentral chi-square distribution. For $r\ge0$, let $\kappa_r$ be the cumulant generating function of $\chi_d^2(rd)$ divided by $d$, and let $J_r(y)$ describe the exponential decay rate of $\Prb\{\chi_d^2(rd)\le yd\}$:
\[
 \kappa_r(t)=-\frac12\log(1-2t)+\frac{rt}{1-2t},\qquad
 {J}_r(y)=\sup_{t<0}\{ty-\kappa_r(t)\},\quad 0<y<1+r.
\]
Set $J_r(y)=0$ for $y\ge1+r$. For $Q\sim\chi_d^2(rd)$, Markov's inequality gives
\[
 {\Prb(Q\le yd)\le e^{-tyd}\E(e^{tQ})
 =\exp\{-d[ty-\kappa_r(t)]\},\qquad t<0.}
\]
Optimizing over $t<0$ gives the upper bound $e^{-dJ_r(y)}$. A matching lower bound on the exponential scale follows by exponential tilting; the Supplement gives the detailed calculation. Thus, for fixed $r\ge0$ and $0<y<1+r$, the noncentral chi-square lower tail satisfies
\[
 {
 \Prb\{\chi_d^2(rd)\le yd\}=\exp\{-dJ_r(y)+o(d)\},\qquad d\to\infty.}
\]

The geometry of $Z_1,Z_2,Z_3$ is described by their normalized Gram matrix, with entries $Z_i^\top Z_j/d$. For a positive-definite $3\times3$ matrix ${V}=({v_{ij}})_{1\le i,j\le3}$, write $w_{ij}({V})={v_{ii}}+{v_{jj}}-2{v_{ij}}$ for the corresponding normalized squared distances. Under the isotropic Gaussian null, the normalized Gram matrix has a Wishart distribution with $d$ degrees of freedom and scale matrix $I_3/d$. Cram\'er's theorem, applied to the iid outer products formed by the coordinate triples $(Z_{1\ell},Z_{2\ell},Z_{3\ell})^\top$ ($\ell=1,\ldots,d$), gives a large-deviation principle with speed $d$ and rate $\frac12\{\tr V-\log{|V|}-3\}$ for $V\succ0$ (see the Supplement).

The next theorem gives the asymptotic behavior of $C_2$ and $C_3$ and the resulting variance orders.

\begin{theorem}
\label{thm:logmotiflaw}
Suppose $d/\log N\to\gamma\in(0,\infty)$ under the isotropic Gaussian null.  Then
\begin{equation}
 {\frac1N\E_0C_2\to0},
 \qquad \frac{C_2}{N}\xrightarrow{L^2}0,
 \label{eq:mainmutuallaw}
\end{equation}
\begin{equation}
 {\E_0C_3}=N^{1+{\beta_\gamma}+o(1)},
 \label{eq:maincollisionlaw}
\end{equation}
where $\beta_\gamma>0$ is given by
\begin{equation}
 {\beta_\gamma}={2-\gamma}\inf_{V\succ0}\bigl\{
 {\tfrac12\{\tr V-\log{|V|}-3\}}\colon
 {w_{21},w_{31}\le w_{23},\min_{j=2,3}J_{v_{jj}}(w_{j1})\ge\gamma^{-1}}\bigr\}.
\label{eq:degreevariationalmain}
\end{equation}
Moreover,
\begin{equation}
 \frac{C_3}{\E_0C_3}\xrightarrow{L^2}1
 \quad\hbox{if }{\frac{(1+\tau_d)^2}{\E_0C_3}}\to0,
\label{eq:mainmotifconcentration}
\end{equation}
where $\tau_d$ is the kissing number in $\mathbb R^d$. Thus, for balanced samples, the permutation variance of $R_w$ is of order $N$, whereas that of $R_d$ is of order ${\E_0C_3}\gg N$, throughout the dimension range in Theorem~\ref{thm:logscalepower}.
\end{theorem}

The conclusion for $C_2$ follows from Theorem~\ref{thm:logendpointgeometry} and exchangeability. The key is that the two endpoints are exchangeable on the mutual-neighbor event, whereas an observation and its selected nearest neighbor have different limiting normalized squared norms. These distinct limits force the probability of mutual pairing to vanish.

Vertices with large in-degrees drive the superlinear growth of ${\E_0}C_3$ and the resulting increase in $\Var_P(R_d)$. To determine this growth rate, we consider the probability that $Z_2$ and $Z_3$ share $Z_1$ as their nearest neighbor. The distance condition $w_{23}\ge\max(w_{21},w_{31})$ prevents either observation from choosing the other in preference to $Z_1$. The two rate constraints arise from requiring all $N-3$ remaining observations to be farther from $Z_2$ and $Z_3$ than $Z_1$ is. The objective function in \eqref{eq:degreevariationalmain} is the large-deviation rate function of the three observations' normalized Gram matrix. Among configurations satisfying these constraints, those with the smallest rate determine the growth of $\E_0C_3$. For each candidate common neighbor, there are order $N^2$ possible ordered pairs of observations. This count gives the term $2$ in \eqref{eq:degreevariationalmain} for the growth exponent of $\E_0C_3/N$. The factor $\gamma$ converts the large-deviation rate per dimension into a decay exponent in $N$, using $d\sim\gamma\log N$; this decay exponent is subtracted from $2$.

The positivity of $\beta_\gamma$ can be seen from a particular feasible matrix. With $r_\gamma=1-s_\gamma+s_\gamma^2$, take
\[
 {V_\gamma=
 \begin{pmatrix}
 r_\gamma&1-s_\gamma&1-s_\gamma\\
 1-s_\gamma&1&(1-s_\gamma)^2/r_\gamma\\
 1-s_\gamma&(1-s_\gamma)^2/r_\gamma&1
 \end{pmatrix}.}
\]
This positive-definite matrix satisfies the constraints in \eqref{eq:degreevariationalmain} and has determinant $s_\gamma^2/r_\gamma$. Substituting it into the rate function gives
\[
 {\tfrac12\{\tr V_\gamma-\log|V_\gamma|-3\}
 =\frac2\gamma-\frac12\{r_\gamma-1-\log r_\gamma\},\qquad
 \beta_\gamma\ge\frac\gamma2\{r_\gamma-1-\log r_\gamma\}>0,}
\]
because $0<r_\gamma<1$.

To control the fluctuations of $C_3$ around its expectation, we bound its change when one observation is replaced. The Efron--Stein inequality then gives
\[
 {\E_0\left|\frac{C_3}{\E_0C_3}-1\right|^2
 \le \frac{C(1+\tau_d)^2}{\E_0C_3}
 \left(1+\frac{N}{\E_0C_3}\right).}
\]
Since $\E_0C_3/N\to\infty$, this proves the relative concentration in \eqref{eq:mainmotifconcentration} under its stated condition. Condition \eqref{eq:mainmotifconcentration} holds automatically in the dimension range of Theorems~\ref{thm:loglocationcomparator} and~\ref{thm:logscalepower}, where $d=(\gamma+o(1))\log N$. The kissing-number bound is $\tau_d\le2^{(0.401413546\ldots+o(1))d}$ as $d\to\infty$ \citep{kabatiansky1978}. Together with \eqref{eq:maincollisionlaw}, this gives $(1+\tau_d)^2/\E_0C_3\le N^{-3/4-\beta_\gamma+o(1)}\to0$. The proof of Theorem~\ref{thm:logmotiflaw} is given in the Supplement.

In the dimension range of Theorems~\ref{thm:loglocationcomparator} and~\ref{thm:logscalepower}, Theorem~\ref{thm:logmotiflaw} gives $C_2/N\to0$ and $C_3/N\to\infty$, while the kissing-number bound gives $C_3/N^2\to0$. Substituting into \eqref{eq:balancedcomponentvariances} gives, in probability,
\[
 {\frac{\Var_P(R_w)}{N}\to\frac1{16},
 \qquad \frac{\Var_P(R_d)}{C_3}\to\frac14.}
\]

The difference $d-\E_0\|Z_{i_*}\|^2$ in Theorem~\ref{thm:logendpointgeometry} controls the first-order scale response of the mean of $R_d$, while the shared-neighbor count in Theorem~\ref{thm:logmotiflaw} determines its fluctuation scale. Their combination gives the local scale for detection in Section~\ref{sec:logscalepower}.

\begin{figure}[!t]
\centering
\includegraphics[width=0.80\textwidth]{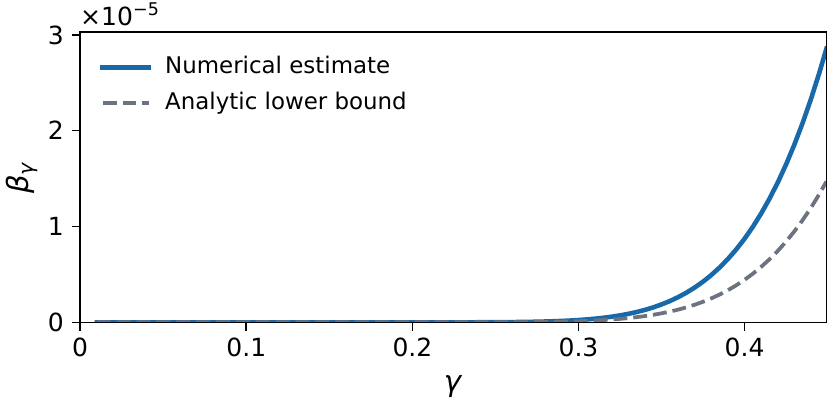}
\caption{Growth exponent of $\E_0C_3/N$ for $d\mathrel{\sim}\gamma\log N$. The solid curve shows numerical estimates of $\beta_\gamma$ and the dashed curve shows the analytic lower bound $\frac{\gamma}{2}\{s_\gamma^2-s_\gamma-\log(1-s_\gamma+s_\gamma^2)\}$, with $s_\gamma$ defined in \eqref{eq:gammasaddle}.}
\label{fig:hubphase}
\end{figure}

Figure~\ref{fig:hubphase} shows numerical estimates of the exponent $\beta_\gamma$ and its analytic lower bound in the dimension range of Theorem~\ref{thm:logscalepower}. Substituting $r_\gamma=1-s_\gamma+s_\gamma^2$ into the lower bound above gives $\frac{\gamma}{2}\{s_\gamma^2-s_\gamma-\log(1-s_\gamma+s_\gamma^2)\}$. In this range, the corresponding local scale is $N^{-1/2+\beta_\gamma/2+o(1)}/d$.

\subsection{\texorpdfstring{{Permutation normal limit and critical value}}{Permutation normal limit and critical value}}
\label{sec:permverification}

{{We use the joint normal limit of} \citet[Theorem~3]{liu2022} to obtain the asymptotic permutation critical value. For the directed 1-NN graph with $m=n=N/2$, the theorem's graph conditions follow from ${d_{\max}}=\max_i{d_i}=o(N^{1/6})$ and $\Var_P(R_d)\ge cN$ for a constant $c>0$ (see the Supplement). Although the approach of \citet{zhu2024} suggests that the degree condition could be substantially relaxed, it already suffices throughout the dimension range of our power theorems.}

In the logarithmic dimension range of Theorems~\ref{thm:loglocationcomparator} and~\ref{thm:logscalepower}, $d/\log N\to\gamma>0$ with $(0.401413546\ldots\log2)\gamma<1/8$. Substituting $d=(\gamma+o(1))\log N$ into the kissing-number bound gives
\[
 {d_{\max}}\le1+\tau_d{{}\le N^{(0.401413546\ldots\log2)\gamma+o(1)}}=o(N^{1/8}).
\]
This satisfies the degree condition for the joint normal limit. For the variance condition, Theorem~\ref{thm:logmotiflaw} and \eqref{eq:balancedcomponentvariances} give $\Var_P(R_d)/N\to\infty$ in probability under the Gaussian null. The covariance comparisons in the proofs of the two power theorems give $\Var_P(R_d)\ge cN$ with probability tending to one under their local alternatives. In one dimension, the degrees are bounded and the variance calculation in the proof of Theorem~\ref{thm:fixedendtoend} gives $\Var_P(R_d)/N\to1/8$ under both the null and its local alternatives.

Thus, under the null and the local alternatives considered in Theorems~\ref{thm:loglocationcomparator}--\ref{thm:fixedendtoend}, the permutation law of $(Z_w,Z_d)^{\top}$ converges in probability to $\N_2(0,I_2)$. Equivalently,
\begin{equation}
\Sigma_P^{-1/2}(R-\mu_P)\Rightarrow\N_2(0,I_2)
\label{eq:conditionalcalibration}
\end{equation}
under the permutation distribution. Since $S=Z_w^2+Z_d^2$, its permutation distribution converges to $\chi_2^2$. The corresponding permutation critical value satisfies $c_{\alpha,N}(G)\to\chi^2_{2,1-\alpha}$ in probability, where $\chi^2_{2,1-\alpha}$ is the $(1-\alpha)$ quantile of the central $\chi^2_2$ distribution.

\section{Power in logarithmic and fixed dimension}
\label{sec:dichotomy}

Consider balanced samples with $F_X=\N(0,I_d)$ and $F_Y=\N(\mu,\sigma^2I_d)$. We first establish a normal approximation under the alternative and then derive local-power limits for location and scale changes when $d$ is proportional to $\log N$. We compare these results with fixed-dimensional connection probabilities and a one-dimensional power limit.

\subsection{\texorpdfstring{Normal approximation {and local power}}{Normal approximation and local power}}
\label{sec:alternativenormal}

{Section~\ref{sec:permverification} gives the asymptotic permutation critical value. To calculate power, we also need the distribution of $S$ under the alternative. The normal approximation below uses $\Cov_a(R)$, whereas $S$ is standardized by $\Sigma_P=\Cov_P(R)$; we therefore compare these matrices after normalization.}

The normal approximation for $R$ and the bound on $\Sigma_P-\E_a\Sigma_P$ both use bounds on changes caused by replacing one observation. We first express these changes in terms of the edges of the two graphs. For two directed graphs $G$ and $G'$ on the same vertex set, $|G\triangle G'|$ counts the directed edges present in exactly one of the two graphs. Write ${d_{\max}}(G)=\max_i{d_i}(G)$ to make the graph explicit. Write $C_2',C_3'$ for the corresponding edge-pair counts in $G'$, and $\|\cdot\|_{\rm op}$ for the operator norm.

\begin{lemma}
\label{lem:edgechanges}
If $G$ and $G'$ are directed graphs without self-loops with the same number of directed edges, and the same sample labels are assigned to their vertices, then
\[
|R_r(G)-R_r(G')|\le {|G\triangle G'|},\quad r=1,2,
\]
and
\[
\begin{gathered}
|C_2-C_2'|\le2{|G\triangle G'|},\qquad
|C_3-C_3'|\le4{|G\triangle G'|}({\max\{{d_{\max}}(G),{d_{\max}}(G')\}}+1),\\
\|\Sigma_P(G)-\Sigma_P(G')\|_{\rm op}\le C_\rho {|G\triangle G'|}({\max\{{d_{\max}}(G),{d_{\max}}(G')\}}+1),
\end{gathered}
\]
whenever $m/N\in[\rho,1-\rho]$ for a fixed $\rho\in(0,1/2)$, where $C_\rho$ depends only on $\rho$.  If $G$ and $G'$ are directed 1-NN graphs obtained by replacing one observation and ${d_{\max}^{\rm in}}$ is the larger maximum in-degree, then
\begin{equation}
{|G\triangle G'|}\le2+4{d_{\max}^{\rm in}}.
\label{eq:editknn}
\end{equation}
\end{lemma}

Each changed edge alters a within-sample count by at most one and the ordered reverse-pair count $C_2$ by at most two. The bound for $C_3$ follows from \eqref{eq:c3degree}, and the covariance bound follows from the exact permutation moments. For \eqref{eq:editknn}, the replaced observation can change its own outgoing edge and the outgoing edges of observations that choose it before or after replacement. Each such change contributes at most two edges to ${G\triangle G'}$. The Supplement gives the full proof of Lemma~\ref{lem:edgechanges}.

For Euclidean directed 1-NN graphs, the kissing-number bound controls both the degrees and the number of edges changed by replacing one observation. Write $d=d_N$ to indicate that dimension may depend on $N$. For a symmetric matrix $A$, let $\lambda_{\min}(A)$ denote its smallest eigenvalue.

\begin{theorem}
\label{thm:alternativepower}
Let $N=2n$, and construct the directed 1-NN graph from independent samples of size $n$ from $\N(0,I_{d_N})$ and $\N(\mu_N,\sigma_N^2I_{d_N})$, where $\sigma_N>0$. Suppose
\[
 {{(1+\tau_{d_N})^3/\sqrt N\to0},\qquad \lambda_{\min}\{\Cov_a(R)\}\ge cN,}
\]
for a constant $c>0$. Then
\[
 {\Cov_a(R)^{-1/2}(R-\E_aR)\Rightarrow\N_2(0,I_2),
 \qquad \|\Sigma_P-\E_a\Sigma_P\|_{\rm op}=o_p(N).}
\]
All stochastic limits are under the alternative sampling distribution. If, in addition, the permutation limit \eqref{eq:conditionalcalibration} holds in probability and
\[
 {\|\Cov_a(R)^{-1/2}\Sigma_P\Cov_a(R)^{-1/2}-I_2\|_{\rm op}\xrightarrow{p}0,}
 \qquad {\Cov_a(R)^{-1/2}}(\E_aR-\mu_P)\xrightarrow{p}\Delta,
\]
for a fixed vector $\Delta\in\mathbb R^2$, then, for $0<\alpha<1$,
\begin{equation}
\Prb_a\{S>c_{\alpha,N}(G)\}\mathrel{\to}
\Prb\{\chi_2^2(\|\Delta\|^2)>\chi^2_{2,1-\alpha}\}.
\label{eq:channelpower}
\end{equation}
Here $\chi_2^2(\lambda)$ denotes a noncentral chi-square distribution with two degrees of freedom and noncentrality parameter $\lambda$.
\end{theorem}

The normal approximation uses \citet[Theorem~2.2]{chatterjee2008}, together with an interaction-graph bound that accounts for the two sample distributions. Lemma~\ref{lem:edgechanges} and the Efron--Stein inequality control the deviation of $\Sigma_P$ from $\E_a\Sigma_P$. The covariance and mean convergence in the theorem then give the power limit by Slutsky's theorem and convergence of the permutation critical value. The proof of Theorem~\ref{thm:alternativepower} is given in the Supplement. For the location and scale alternatives below, we verify the covariance comparison in Theorem~\ref{thm:alternativepower} and calculate the limiting standardized mean $\Delta$.

The dimension range in the next two subsections is expressed using the constant $\kappa_{\rm KL}$ in the kissing-number bound. Put $s=\sqrt3/2$, $a=(1+s)/(2s)$, and $b=(1-s)/(2s)$, and define
\[
 \kappa_{\rm KL}
 =\{a\log_2a-b\log_2b\}\log2
 =0.401413546\ldots\log2.
\]
The kissing-number bound $\tau_d\le\exp\{(\kappa_{\rm KL}+o(1))d\}$ follows from the spherical-code bound of \citet{kabatiansky1978} at minimum angle $\pi/3$. This bound implies the degree condition in Theorem~\ref{thm:alternativepower}, $(1+\tau_{d_N})^3/\sqrt N\to0$, whenever $d_N\le c\log N$ eventually for some $c<(6\kappa_{\rm KL})^{-1}$. Indeed,
\[
 {\frac{(1+\tau_{d_N})^3}{\sqrt N}
 \le N^{3\kappa_{\rm KL}c-1/2+o(1)}\mathrel{\to}0.}
\]
For the covariance comparisons in the power proofs below, we use the narrower range $\gamma<(8\kappa_{\rm KL})^{-1}$.

\subsection{Logarithmic Gaussian location power}

We first consider location alternatives when $d/\log N\to\gamma>0$. For balanced Gaussian samples, the symmetry in \eqref{eq:exactlocationchannels} gives $\E_aR_d=0$ at every sample size. We now calculate the common change in the two connection probabilities and compare the resulting change in $R_w$ with its standard deviation.

Write $p_{r,N}$ for the within-sample connection probability $p_r$, $r=1,2$, defined in Section~\ref{sec:exact}, to make its dependence on sample size explicit. For $b_N>0$, we write $a_N=O(b_N)$ if $|a_N|\le Cb_N$ for all sufficiently large $N$, with a constant $C$ independent of $N$.

\begin{theorem}
\label{thm:loglocationcomparator}
Let $N=2n$, let $d/\log N\allowbreak\to\gamma\in (0,(8\kappa_{\rm KL})^{-1})$, and consider $F_X=\N(0,I_d)$ and $F_Y=\N(\mu_N,I_d)$.  Write $p_{0,N}=(n-1)/(2n-1)$ for the null within-sample connection probability.  For every fixed $C>0$, uniformly over $\|\mu_N\|^2\mathrel{\le}CN^{-1/2}$,
\begin{equation}
 p_{1,N}=p_{2,N}
 =p_{0,N}+\frac{1-s_\gamma}{8}\|\mu_N\|^2
 +o(\|\mu_N\|^2)+O(\|\mu_N\|^4).
\label{eq:loglocationcomparatormean}
\end{equation}
If ${N^{1/4}\|\mu_N\|\to h\ge0}$, the two components satisfy
\[
 {(Z_w,Z_d)^{\top}\Rightarrow
 \N_2\!\left\{\left(\frac{h^2(1-s_\gamma)}4,0\right)^{\top},I_2\right\}.}
\]
The level-$\alpha$ upper-tail permutation test based on $Z_d$ therefore has limiting power $\alpha$, for $0<\alpha<1$. The test based on $S$ using the exact permutation covariance and critical value satisfies
\begin{equation}
 \Prb_{\mu_N}\{S>c_{\alpha,N}(G)\}
 \mathrel{\to}
 \Prb\!\left\{
 \chi_2^2\!\left(\frac{{h^4}(1-s_\gamma)^2}{16}\right)
 >\chi^2_{2,1-\alpha}\right\}.
\label{eq:loglocationcomparatorpower}
\end{equation}
\end{theorem}

At $N^{1/4}\|\mu_N\|\to h$, this expansion gives a mean shift of $\sqrt N h^2(1-s_\gamma)/16+o(\sqrt N)$ in $R_w$. Under the null, Theorem~\ref{thm:logmotiflaw} and \eqref{eq:balancedcomponentvariances} give $\Var_P(R_w)/N\to1/16$ in probability. Location shifts of order $N^{-1/4}$ therefore produce a nonvanishing standardized mean shift in $Z_w$, giving the weighted component nontrivial local power.

We next explain the main steps in the proof of Theorem~\ref{thm:loglocationcomparator}, beginning with the quadratic expansion of the connection probabilities. Since $p_{1,N}=p_{2,N}$, it is enough to consider an $X$ observation. Write $\mu=tv$ with $\|v\|=1$, and let $p_1^{\rm loc}(t)$ denote this probability. For a source $X_1$, let $W=\min_{2\le i\le n}\|X_i-X_1\|^2$ and define $F_t(w\mid x)=\Prb\{\|Y-x\|^2\le w\}$ for an independent $Y\sim\N(tv,I_d)$. The nearest neighbor of $X_1$ belongs to sample $X$ exactly when all $n$ observations from sample $Y$ lie outside this distance ball. Conditional on $X_1$ and $W$, these $n$ events are independent, so
\[
 {p_1^{\rm loc}(t)=\E\{[1-F_t(W\mid X_1)]^n\}.}
\]
Put $U=F_0(W\mid X_1)$. Conditional on $X_1$, the variable $U$ has the $\operatorname{Beta}(1,n-1)$ distribution. Write $\dot F_0$ and $\ddot F_0$ for the first two derivatives with respect to $t$ at zero. Twice differentiating the connection probability therefore gives
\[
 {
 \{p_1^{\rm loc}\}''(0)
 =\E[n(n-1)(1-U)^{n-2}\dot F_0(W\mid X_1)^2]-\E[n(1-U)^{n-1}\ddot F_0(W\mid X_1)].
 }
\]
The noncentral chi-square tail expansion and the beta law give the limits of the two terms:
\[
 {\begin{aligned}
 \E[n(n-1)(1-U)^{n-2}\dot F_0(W\mid X_1)^2]
 &\mathrel{\to} \frac{(1-s_\gamma)^2}{4},\\
 \E[n(1-U)^{n-1}\ddot F_0(W\mid X_1)]
 &\mathrel{\to} \frac{s_\gamma-1+(1-s_\gamma)^2}{4}.
 \end{aligned}}
\]
Subtracting gives $\{p_1^{\rm loc}\}''(0)\to(1-s_\gamma)/4$. Reflection symmetry makes the odd derivatives at zero vanish, and a uniform fourth-derivative bound controls the Taylor remainder. This gives \eqref{eq:loglocationcomparatormean}.

To verify the covariance conditions in Theorem~\ref{thm:alternativepower}, we compare both $\Cov_a(R)$ and $\E_a\Sigma_P$ with $\E_0\Sigma_P$. Write $\Cov_0(R)$ for the covariance of $R$ under the null sampling distribution. Under the null, the observations are iid and the graph rule is equivariant under vertex relabeling, so every fixed allocation of $m$ sample-$X$ labels and $n$ sample-$Y$ labels gives the same sampling distribution of $R$. Since $\E_P$ averages over these allocations for each observed graph, the null first and second moments satisfy
\[
 {\begin{aligned}
 \E_0R&=\E_0(\E_PR)=\E_0\mu_P,\\
 \E_0(RR^\top)&=\E_0\{\E_P(RR^\top)\}
 =\E_0\{\Sigma_P+\mu_P\mu_P^\top\}.
 \end{aligned}}
\]
For fixed sample sizes and $|G|=N$, $\mu_P$ is deterministic. Subtracting $\mu_P\mu_P^\top$ from the second identity gives $\Cov_0(R)=\E_0\Sigma_P$. Theorem~\ref{thm:logmotiflaw} and the exact variance formulas \eqref{eq:balancedcomponentvariances} show that $\E_0\Sigma_P$ has eigenvalues asymptotic to $N/8$ and $\E_0C_3/8$ in the unit sum and difference directions. Since $\E_0C_3/N\to\infty$, the smaller eigenvalue of $\E_0\Sigma_P$ is $N/8+o(N)$.

Reflecting all observations about the origin changes the location shift from $tv$ to $-tv$ while leaving $R$ and $\Sigma_P$ unchanged. Consequently, both the sampling covariance and the expected permutation covariance are even functions of $t$, with first derivatives zero at $t=0$. Gaussian score differentiation and the graph-degree bounds show that every entry of $\Cov_a(R)$ and $\E_a\Sigma_P$, viewed as a function of $t$ along $\mu=tv$, has second derivative bounded in absolute value by $CN(1+\tau_d)^4\log(eN)$, uniformly for $|t|\le\|\mu_N\|$. Together with Taylor's theorem, these bounds give
\[
 {\begin{aligned}
 &\frac{\|\Cov_a(R)-\E_0\Sigma_P\|_{\rm op}
       +\|\E_a\Sigma_P-\E_0\Sigma_P\|_{\rm op}}{N}\\
 &\qquad\le C\|\mu_N\|^2(1+\tau_d)^4\log(eN)
 \le N^{4\kappa_{\rm KL}\gamma-1/2+o(1)}\mathrel{\to}0.
 \end{aligned}}
\]
The last step uses $\|\mu_N\|^2=O(N^{-1/2})$ and $\gamma<(8\kappa_{\rm KL})^{-1}$. Thus $\lambda_{\min}\{\Cov_a(R)\}/N\to1/8$, which verifies the variance lower bound in Theorem~\ref{thm:alternativepower}. Lemma~\ref{lem:edgechanges} and the Efron--Stein inequality also give $\|\Sigma_P-\E_a\Sigma_P\|_{\rm op}=O_p\{\sqrt N(1+\tau_d)^2\}=o_p(N)$. Combining this fluctuation bound with the two covariance comparisons yields
\[
 {\|\Cov_a(R)^{-1/2}\Sigma_P\Cov_a(R)^{-1/2}-I_2\|_{\rm op}
 \le\frac{\|\Sigma_P-\Cov_a(R)\|_{\rm op}}
 {\lambda_{\min}\{\Cov_a(R)\}}
 =o_p(1).}
\]

It remains to calculate the limiting standardized mean $\Delta$ in the original $(R_1,R_2)$ coordinates. By balanced location symmetry, the mean shift lies in the sum direction, which is an eigendirection of $\Cov_a(R)$ with eigenvalue $2\Var_a(R_w)=N/8+o(N)$. Applying $\Cov_a(R)^{-1/2}$ therefore divides the mean shift by $\sqrt{N/8}\{1+o(1)\}$. Substituting the mean expansion gives
\[
 {\Cov_a(R)^{-1/2}(\E_aR-\mu_P)
 =\frac{\sqrt N\,h^2(1-s_\gamma)/16}{\sqrt{N/8}}(1,1)^\top+o(1)
 \mathrel{=}\frac{h^2(1-s_\gamma)}{4\sqrt2}(1,1)^\top{+o(1)}.}
\]
The orthogonal change to sum and difference coordinates takes the limiting vector $\Delta$ to $(h^2(1-s_\gamma)/4,0)^{\top}$, the limiting drifts of $(Z_w,Z_d)^{\top}$. The normal approximation under the alternative and the permutation limit then give \eqref{eq:loglocationcomparatorpower}. The detailed proof of the theorem is given in the Supplement.

\subsection{Logarithmic Gaussian scale power}
\label{sec:logscalepower}

Write $p_r^{\rm sc}(\varepsilon)$ for the within-sample connection probability in sample $r$ when $\sigma=1+\varepsilon$ and $\mu=0$. Under the balanced null, let ${Z_1},{Z_2},\ldots,{Z_N}\overset{\mathrm{iid}}{\sim}\N(0,I_d)$, where $N=2n$, and let ${i_*}=\arg\min_{2\le j\le N}\|Z_j-{Z_1}\|$. For the $\N\{0,(1+\varepsilon)^2I_d\}$ density, the log-density derivative at $\varepsilon=0$ is $\|y\|^2-d$. Differentiating $p_1^{\rm sc}$ therefore gives the null expectation of the within-sample connection indicator for an $X$ observation multiplied by $\sum_{j=1}^n(\|Y_j\|^2-d)$.

For this null expectation, keep the label of $Z_1$ fixed as $X$ and independently permute the remaining labels uniformly, preserving the two sample sizes. By exchangeability, the expectation is unchanged when both the connection indicator and the scale score are evaluated with these permuted labels. Conditional on the pooled observations, the nearest neighbor of $Z_1$ has label $X$ with probability $(n-1)/(N-1)$. Conditional also on the nearest neighbor of $Z_1$ having label $X$, the $n$ observations assigned to sample $Y$ form a uniform subset of the other $N-2$ observations. Since
\[
 {\E_0\!\Bigl[\sum_{2\le i\le N,\ i\ne i_*}(\|Z_i\|^2-d)\Bigr]
 =d-\E_0\|Z_{i_*}\|^2,}
\]
the label average gives $\{p_1^{\rm sc}\}'(0)$ by multiplying this expectation by $(n-1)n/\{(N-1)(N-2)\}=n/\{2(2n-1)\}$. Interchanging the samples and rescaling all observations by $(1+\varepsilon)^{-1}$ changes the scale parameter to $-\varepsilon/(1+\varepsilon)$, so the two derivatives at zero have opposite signs. Thus, for every finite $(n,d)$,
\begin{equation}
 \{p_1^{\rm sc}\}'(0)=-\{p_2^{\rm sc}\}'(0)
 =\frac{n}{2(2n-1)}\{d-\E_0\|Z_{i_*}\|^2\}.
\label{eq:scalenearestnorm}
\end{equation}
The Supplement gives the full derivation of \eqref{eq:scalenearestnorm}.

Since $R_d=R_1-R_2$ and $m=n$, we have $\E_aR_d=n\{p_1^{\rm sc}(\varepsilon)-p_2^{\rm sc}(\varepsilon)\}$. Differentiating gives
\[
 {\left.\frac{\mathrm d}{\mathrm d\varepsilon}\E_aR_d\right|_{\varepsilon=0}
 =n\bigl\{\{p_1^{\rm sc}\}'(0)-\{p_2^{\rm sc}\}'(0)\bigr\}
 =N\{p_1^{\rm sc}\}'(0).}
\]
Theorem~\ref{thm:logendpointgeometry} gives $\{p_1^{\rm sc}\}'(0)/d\to s_\gamma(1-s_\gamma)/4>0$ when $d/\log N\to\gamma>0$. Thus the linear term in the scale response of $\E_{a} R_d$ is of order $\varepsilon_NNd$.

Section~\ref{sec:permutationmoments} shows that the permutation standard deviation of $R_d$ grows faster than $\sqrt N$ in the logarithmic dimension range of our power results. Comparing the linear mean term with $\sqrt N$ would suggest the scale $1/(\sqrt N d)$. The next proposition rules out asymptotic power above the test level at this scale. More generally, it rules out such power whenever $|\varepsilon_N|=o\{(Nd)^{-1/2}\}$.

\begin{proposition}
\phantomsection
\label{prop:compatibility}
Let $N=2n\to\infty$ and let $d=d_N\ge1$ be an integer. Consider independent samples
\[
 {X_1,\ldots,X_n\overset{\mathrm{iid}}{\sim}\N(0,I_d),\qquad
 Y_1,\ldots,Y_n\overset{\mathrm{iid}}{\sim}\N\{0,(1+\varepsilon_N)^2I_d\},}
\]
where $1+\varepsilon_N>0$. Write $\Prb_{\varepsilon,N}$ for their joint distribution and $\Prb_{0,N}$ for their joint distribution under the null hypothesis. If $Nd\varepsilon_N^2\to0$, then
\[
D_{\rm KL}(\Prb_{\varepsilon,N}\Vert\Prb_{0,N})\to0,\qquad
\|\Prb_{\varepsilon,N}-\Prb_{0,N}\|_{\rm TV}\to0.
\]
Consequently{,} every sequence of tests with asymptotic size at most $\alpha$ has power with limit superior at most $\alpha$, for $0<\alpha<1$. Thus
\[
 |\varepsilon_N|\sqrt{Nd}\not\to0
\]
is necessary for asymptotic power above the test level.
\end{proposition}

The power bound in Proposition~\ref{prop:compatibility} follows by comparing the joint distributions under the null and alternative. Only the distribution of the second sample changes, so the joint Kullback--Leibler divergence is the sum of the contributions from its $n=N/2$ independent observations. Since $Nd\varepsilon_N^2\to0$ implies $\varepsilon_N\to0$, the Gaussian divergence satisfies
\[
 {\begin{aligned}
 D_{\rm KL}(\Prb_{\varepsilon,N}\Vert\Prb_{0,N})
 &=\frac{nd}{2}\left[(1+\varepsilon_N)^2-1-\log\{(1+\varepsilon_N)^2\}\right]\\
 &=nd\varepsilon_N^2\{1+O(|\varepsilon_N|)\}\mathrel{\to}0.
 \end{aligned}}
\]
Pinsker's inequality then gives
\[
 {\|\Prb_{\varepsilon,N}-\Prb_{0,N}\|_{\rm TV}
 \le\sqrt{\tfrac12D_{\rm KL}(\Prb_{\varepsilon,N}\Vert\Prb_{0,N})}\mathrel{\to}0.}
\]
For any test, including a randomized test, the absolute difference between its null and alternative rejection probabilities is bounded by this total-variation distance. Taking limit superiors proves the power bound.

For the directed 1-NN test, comparing the linear mean term of order $\varepsilon_NNd$ with the permutation standard deviation of order $\sqrt{\E_0C_3}$ gives the normalization used in the following theorem.

\begin{theorem}
\label{thm:logscalepower}
Let $N=2n$, let $d/\log N\allowbreak\to\gamma\in (0,(8\kappa_{\rm KL})^{-1})$, and consider $F_X=\N(0,I_d)$ and $F_Y=\N\{0,(1+\varepsilon_N)^2I_d\}$. If
\begin{equation}
 \frac{\varepsilon_NNd}{\sqrt{\E_0C_3}}\mathrel{\to}\ell,
 \label{eq:hubscalesequence}
\end{equation}
for some $\ell\in\mathbb R$, then the two components satisfy
\[
 {(Z_w,Z_d)^{\top}\Rightarrow
 \N_2\!\left\{\left(0,\frac{\ell s_\gamma(1-s_\gamma)}2\right)^{\top},I_2\right\}.}
\]
The level-$\alpha$ upper-tail permutation test based on $Z_w$ therefore has limiting power $\alpha$, for $0<\alpha<1$. The test based on $S$ using the exact permutation covariance and critical value satisfies
\begin{equation}
 \Prb_{\varepsilon,N}\{S>c_{\alpha,N}(G)\}
 \mathrel{\to}
 \Prb\!\left\{
 \chi_2^2\!\left(
 \frac{\ell^2s_\gamma^2(1-s_\gamma)^2}{4}
 \right)>\chi^2_{2,1-\alpha}\right\}.
 \label{eq:balancedscalepower}
\end{equation}
\end{theorem}

We next explain the main steps in the proof of Theorem~\ref{thm:logscalepower}, beginning with the mean shifts of the two components. For the difference count, equation~\eqref{eq:scalenearestnorm} and Theorem~\ref{thm:logendpointgeometry} determine the linear term in $\E_aR_d$. Gaussian score identities express the second derivatives of the count means in terms of their joint cumulants with the scale score. The graph-degree bounds control how many edges can change when observations are replaced. Together with the Gaussian score moment bounds, they control the cumulants and hence the Taylor remainders for both components. On the local sequence \eqref{eq:hubscalesequence}, these bounds give
\[
 {\E_aR_d=\frac{\varepsilon_NNd\,s_\gamma(1-s_\gamma)}4
       +o\!\left(\sqrt{\E_0C_3}\right).}
\]

For $R_w$, the first derivative of $\E_aR_w$ with respect to $\varepsilon$ is zero at $\varepsilon=0$. Taylor's theorem and the same degree bounds give, on this local sequence,
\[
 {|\E_aR_w-\E_PR_w|
 \le \frac{C(1+\tau_d)^3\E_0C_3}{Nd}
 \le \frac{C(1+\tau_d)^4}{d}=o(\sqrt N).}
\]
Here $\E_PR_w=\E_0R_w$ for every directed 1-NN graph with balanced samples. The second inequality uses $\E_0C_3\le 2N(1+\tau_d)$, and the final limit follows from $1+\tau_d=o(N^{1/8})$ in the theorem's dimension range. Thus the remaining weighted mean shift is negligible relative to its standard deviation of order $\sqrt N$.

To verify the covariance conditions in Theorem~\ref{thm:alternativepower}, we again compare $\Cov_a(R)$ and $\E_a\Sigma_P$ with $\E_0\Sigma_P$. We parameterize scale by $\eta=\varepsilon\sqrt d$, with $\eta_N=\varepsilon_N\sqrt d$. Under this parameterization, the fourth moment of each scaled observation's score is bounded uniformly in $d$ along the local path. Gaussian score differentiation and the graph-degree bounds show that every entry of $\Cov_a(R)$ and $\E_a\Sigma_P$, viewed as a function of $\eta$ with $\varepsilon=\eta/\sqrt d$, has first derivative bounded in absolute value by $CN(1+\tau_d)^3\log(eN)$, uniformly for $|\eta|\le|\eta_N|$. Integrating from the null to $\eta_N$ gives
\[
 {\begin{aligned}
 &\frac{\|\Cov_a(R)-\E_0\Sigma_P\|_{\rm op}
       +\|\E_a\Sigma_P-\E_0\Sigma_P\|_{\rm op}}{N}\\
 &\qquad\le C|\varepsilon_N|\sqrt d(1+\tau_d)^3\log(eN)
 \le\frac{C(1+\tau_d)^{7/2}\log(eN)}{\sqrt{Nd}}
 \mathrel{\to}0.
 \end{aligned}}
\]
The second inequality uses the local scale \eqref{eq:hubscalesequence} and $\E_0C_3\le2N(1+\tau_d)$. The limit follows from $1+\tau_d=o(N^{1/8})$ and $d\asymp\log N$. As in the location calculation, comparison with $\E_0\Sigma_P$ gives $\lambda_{\min}\{\Cov_a(R)\}/N\to1/8$, verifying the variance lower bound in Theorem~\ref{thm:alternativepower}. Lemma~\ref{lem:edgechanges} and the Efron--Stein inequality also give $\|\Sigma_P-\E_a\Sigma_P\|_{\rm op}=o_p(N)$ under the scale alternative. Combining this fluctuation bound with the two covariance comparisons yields $\|\Sigma_P-\Cov_a(R)\|_{\rm op}=o_p(N)$. Together with the variance lower bound, this verifies the normalized covariance condition in Theorem~\ref{thm:alternativepower}. For the permutation variances, these bounds give $\sigma_w^2/N\to1/16$ and $\sigma_d^2/\E_0C_3\to1/4$ in probability under the alternative.

It remains to calculate the limiting standardized mean $\Delta$ in the original $(R_1,R_2)$ coordinates. The sampling variances of the projections of $R$ onto the unit sum and difference directions satisfy $2\Var_a(R_w)\sim N/8$ and $\Var_a(R_d)/2\sim\E_0C_3/8$, respectively, and the covariance between these directions is $\Cov_a(R_w,R_d)=o(\sqrt{N\E_0C_3})$. Thus $\Cov_a(R)^{-1/2}$ asymptotically divides the mean shift in these directions by $\sqrt{N/8}$ and $\sqrt{\E_0C_3/8}$. Substituting the two mean expansions above gives
\[
 {\Cov_a(R)^{-1/2}(\E_aR-\mu_P)
 =\frac{\varepsilon_NNd\,s_\gamma(1-s_\gamma)}{\sqrt{8\E_0C_3}}{\begin{pmatrix}1\\-1\end{pmatrix}}+o(1)
 \mathrel{=}\frac{\ell s_\gamma(1-s_\gamma)}{2\sqrt2}{\begin{pmatrix}1\\-1\end{pmatrix}}{+o(1)}.}
\]
The orthogonal change to sum and difference coordinates takes this limiting vector to $(0,\ell s_\gamma(1-s_\gamma)/2)^{\top}$, the limiting drifts of $(Z_w,Z_d)^{\top}$. Combining the centered normal approximation under the alternative with these shifts gives the joint limit in the theorem. Together with the permutation critical-value limits, this gives limiting power $\alpha$ for the upper-tail test based on $Z_w$ and the power limit \eqref{eq:balancedscalepower} for $S$. The detailed proof of the theorem is given in the Supplement.

\subsection{A fixed-dimensional comparison}

We now consider balanced samples in fixed dimension. The distance to a nearest neighbor tends to zero. Let $f$ and $f_t$ be the two sample densities on $\R^d$, with $t$ indexing the alternative. The limiting within-sample connection probabilities are given by \citet[eqs.~(3.5)--(3.6)]{schilling1986}:
\begin{equation}
 {p_{1,\infty}}(t)=\int\frac{f^2}{f+f_t},\qquad
 {p_{2,\infty}}(t)=\int\frac{f_t^2}{f+f_t}.
 \label{eq:fixedoverlaps}
\end{equation}
The integrands are taken to be zero where $f+f_t=0$. Subtracting the two integrals in \eqref{eq:fixedoverlaps} gives
\[
 {p_{1,\infty}(t)-p_{2,\infty}(t)
 =\int\frac{f^2-f_t^2}{f+f_t}=\int(f-f_t)=0.}
\]
This equality holds for both location and scale alternatives.

To calculate their common response to a scale change, consider the log-scale path parameterized by $t=\log(1+\varepsilon)$:
\begin{equation}
 f_t(x)=e^{-dt}f(e^{-t}x),\qquad
 \dot\ell_{\rm sc}(x)=-d-x^{\top}\nabla\log f(x).
 \label{eq:fixedscalepath}
\end{equation}
The score is defined where $f>0$. Write $I_{\rm sc}(f)=\int_{\{f>0\}} \dot\ell_{\rm sc}^2 f$. Define $u_t=(f_t-f)/f$ on $\{f>0\}$. Suppose $0<I_{\rm sc}(f)<\infty$ and, with the following integrals taken over this set,
\begin{equation}
 \int f\{u_t-t\dot\ell_{\rm sc}\}^2=o(t^2),\qquad
 \int f|u_t|^3=o(t^2).
 \label{eq:fixedoverlapregularity}
\end{equation}
The equality above and a Taylor expansion of \eqref{eq:fixedoverlaps} give, for every sufficiently small $t$ for which the displayed integrals are finite,
\begin{equation}
 {p_{1,\infty}}(t)-{p_{2,\infty}}(t)=0,\qquad
 {p_{r,\infty}}(t)=\frac12+\frac{t^2}{8}I_{\rm sc}(f)+o(t^2),\quad r=1,2.
 \label{eq:fixedoverlapexpansion}
\end{equation}
Thus, for every fixed $d$, the two limiting connection probabilities have the same quadratic response to a small scale change.

For scale changes, the finite-sample connection probabilities also satisfy an exact symmetry.  For balanced samples, let $p_{r,N}(t)$ denote the within-sample connection probability in sample $r$ under $(f,f_t)$. Multiplying every observation by the same positive constant leaves the nearest-neighbor graph unchanged.  If data from $(f,f_t)$ are multiplied by $e^{-t}$ and the two sample labels are then interchanged, the two sample distributions become $(f,f_{-t})$. Therefore
\begin{equation}
 p_{1,N}(t)=p_{2,N}(-t)
 \label{eq:fixedfiniteparity}
\end{equation}
for every $N$, $d$, and admissible $t$.  Thus the weighted sum has no first-order response to a scale change.

For a location shift, let $f_t(x)=f(x-tv)$ with $\|v\|=1$. The expansion in \eqref{eq:fixedoverlapexpansion} holds with $\dot\ell_{\rm sc}(x)$ replaced by $-v^\top\nabla\log f(x)$ and $I_{\rm sc}(f)$ by $\int_{\{f>0\}}(v^\top\nabla\log f)^2f$, provided this integral is finite and positive and \eqref{eq:fixedoverlapregularity} holds for the location path and its score. For $f$ the $\N(0,I_d)$ density, this integral equals one, so both limiting connection probabilities increase by $t^2/8+o(t^2)$.

For either family, the quadratic overlap expansion suggests a count mean shift of order $Nt_N^2$. Comparing it with a standard deviation of order $\sqrt N$ gives $t_N\asymp N^{-1/4}$. For balanced Gaussian samples in one dimension, the following theorem gives the local-power limits. \begingroup

\begin{theorem}
\label{thm:fixedendtoend}
Let the two samples be independent, with iid observations within each sample. Let $N=2n$ and construct the directed $1$-NN graph on the line.  Consider either
\[
\begin{array}{lll}
\text{\rm location:}&X_i\sim\N(0,1),&
Y_j\sim\N(N^{-1/4}h,1),\quad h\in\mathbb R,\\[2mm]
\text{\rm scale:}&X_i\sim\N(0,1),&
Y_j\sim\N\{0,(1+aN^{-1/4})^2\},\quad a\in\mathbb R.
\end{array}
\]
Set $\Delta_w=h^2\sqrt{3/80}$ in the location case and $\Delta_w=a^2\sqrt{3/20}$ in the scale case.  Then the two permutation-standardized components satisfy
\begin{equation}
 (Z_w,Z_d)^{\top}\Rightarrow
 \N_2\{(\Delta_w,0)^{\top},I_2\}.
 \label{eq:onedjointchannels}
\end{equation}
In particular, the weighted component has nonzero limiting drift when $\Delta_w>0$, whereas the difference component has zero limiting drift. Under the corresponding alternative, the exact permutation-studentized statistic satisfies
\begin{equation}
 {\Prb}\{S>c_{\alpha,N}(G)\}
 \mathrel{\to}
 \Prb\left[
 \chi_2^2\left(\Delta_w^2\right)
 >\chi^2_{2,1-\alpha}\right].
\label{eq:fixedpowercurve}
\end{equation}
The location and scale noncentralities are $3h^4/80$ and $3a^4/20$, respectively.
\end{theorem}

\endgroup {For location, the $N^{-1/4}$ local scale and the component carrying the drift agree with Theorem~\ref{thm:loglocationcomparator}, although the constants differ. For scale, the nonzero limiting drift is in $Z_w$ in one dimension and in $Z_d$ in Theorem~\ref{thm:logscalepower}, at their respective local scales.}

To prove Theorem~\ref{thm:fixedendtoend}, we first approximate the finite-sample connection probabilities by the density-overlap integrals with error $o(N^{-1/2})$, uniformly over the local alternatives in the theorem. In one dimension, this follows by localizing the exact nearest-distance integral and controlling the Gaussian tails. The quadratic overlap expansions then give the common mean shift with an $o(\sqrt N)$ error in the edge counts, while $\E_aR_d=o(\sqrt N)$. The one-dimensional graph-count limits, covariance comparison, and normal approximations give the stated standardized shifts and power. The detailed proof of the theorem is given in the Supplement.

\section{Numerical illustration}
\label{sec:sims}

We first examine how connection probabilities change under small location shifts and scale changes when dimension is proportional to $\log N$, then compare connection probabilities, component means, and power predictions for both location and scale alternatives at a fixed sample size.

\subsection{\texorpdfstring{Connection probabilities under small location and scale changes}{Connection probabilities under small location and scale changes}}

For balanced Gaussian samples, the leading change in $p_1^{\rm loc}(t)$ is quadratic in the shift magnitude $t=\|\mu\|$, whereas that in $p_1^{\rm sc}(\varepsilon)$ is linear in the relative scale change $\varepsilon$. Theorem~\ref{thm:loglocationcomparator} and \eqref{eq:scalenearestnorm} give the following limits for the coefficient of $t^2$ and the scale derivative multiplied by $4/d$
\[
\begin{aligned}
 {c_{N,d}}&\mathrel{=}{\frac12\{p_1^{\rm loc}\}''(0)
 \to\frac{1-s_\gamma}{8}},\\
 a_{N,d}&=\frac4d\{p_1^{\rm sc}\}'(0)
 =\frac{2n}{2n-1}\left(1-\frac{\E_0\|Z_{i_*}\|^2}{d}\right)
 \mathrel{\to}s_\gamma(1-s_\gamma).
\end{aligned}
\]

In Table~\ref{tab:theoremwindow}, we set $\gamma=0.4$ and take $N$ to be the nearest even integer to $\exp(d/0.4)$; the table reports the resulting value of $d/\log N$. We obtain $s_\gamma$ for the limiting values by solving \eqref{eq:gammasaddle} numerically with $\gamma=d/\log N$. We evaluate $c_{N,d}$ by numerical integration of the exact second-derivative formula used in the proof of Theorem~\ref{thm:loglocationcomparator}. For $a_{N,d}$, conditional on the source observation, we sample the minimum squared distance using the exact noncentral-chi-square distribution of each candidate squared distance. The squared norm of the selected nearest neighbor is replaced by its conditional expectation under the exact angular law. The finite-sample values of $c_{N,d}$ and $a_{N,d}$ move toward their predicted limits over the displayed range, with a smaller relative discrepancy for location.

\begin{table}[!t]
\centering
\small
\caption{Quadratic location coefficient $c_{N,d}$ and scale derivative multiplied by $4/d$ ($a_{N,d}$) for $d/\log N\approx0.4$. Location entries are computed by numerical integration; MCSE denotes Monte Carlo standard error for the scale estimates.}
\label{tab:theoremwindow}
\begin{tabular}{rrrrrrrr}
\toprule
\multicolumn{3}{c}{} & \multicolumn{2}{c}{Location} & \multicolumn{3}{c}{Scale}\\
\cmidrule(lr){4-5}\cmidrule(lr){6-8}
$d$ & $\log N$ & $d/\log N$ & $c_{N,d}$ & $(1-s_\gamma)/8$ & $a_{N,d}$ & MCSE & $s_\gamma(1-s_\gamma)$\\
\midrule
3 & 7.5 & 0.400 & 0.120897 & 0.124163 & 0.04733 & 0.00071 & 0.00665 \\
6 & 15.0 & 0.400 & 0.122939 & 0.124163 & 0.02122 & 0.00054 & 0.00665 \\
12 & 30.0 & 0.400 & 0.123643 & 0.124163 & 0.01236 & 0.00039 & 0.00665 \\
24 & 60.0 & 0.400 & 0.123913 & 0.124163 & 0.00937 & 0.00028 & 0.00665 \\
 \bottomrule
\end{tabular}
\end{table}

\subsection{\texorpdfstring{Connection probabilities and component means}{Connection probabilities and component means}}
\label{sec:finiteconnection}
\suppressfloats[t]

{We examine how finite-sample connection probabilities translate into mean shifts in $Z_w$ and $Z_d$.} We fix the sample sizes at $m=n=100$ ($N=200$) and vary $d$ over 2, 5, 10, 20, 50, and 100. The Gaussian alternatives have $\|\mu\|=2$ for location and $\sigma=1.08$ for scale.

\begin{figure}[!b]
\centering
\includegraphics[width=0.94\textwidth]{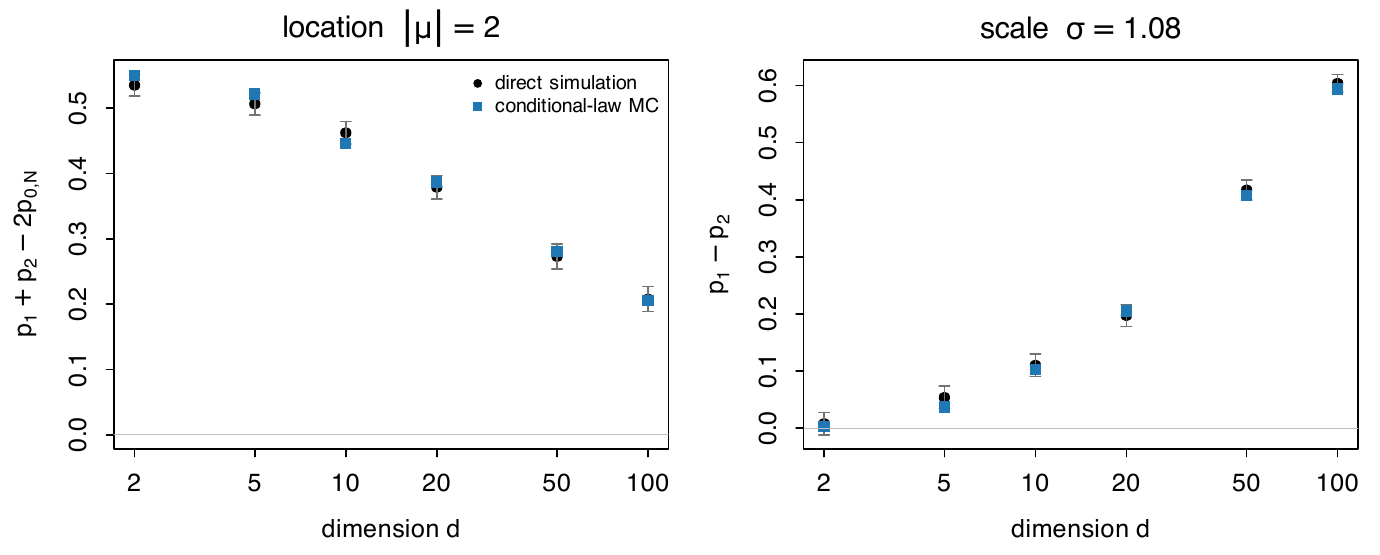}
\caption{Sums and differences of within-sample connection probabilities at $N=200$. Here $p_{0,N}=(n-1)/(2n-1)$ is the null connection probability. Circles show direct simulation; squares show conditional-law Monte Carlo (MC) estimates. Error bars show $\pm1.96$ Monte Carlo standard errors for the direct estimates.}
\label{fig:dichotomy}
\end{figure}

Figure~\ref{fig:dichotomy} shows that $p_1+p_2-2p_{0,N}$ decreases with dimension under the location alternative. Under the scale alternative, $p_1-p_2$ is indistinguishable from zero at Monte Carlo resolution when $d=2$ and then grows steadily across the displayed dimensions.

We estimate the connection probabilities in two ways: by direct simulation of the observations and by conditional-law Monte Carlo (MC). The conditional-law MC calculation first draws observations from each sample. Given an observation, it generates the within-sample and between-sample squared distances independently from the exact scaled noncentral chi-square laws in Section~\ref{sec:gaussian}. It records whether the within-sample minimum is smaller and averages these indicators over repeated distance draws and sampled observations to estimate $p_1$ and $p_2$. The conditional-law Monte Carlo estimates closely follow the direct estimates.

We next use the connection probabilities to predict the means of $Z_w$ and $Z_d$ at $N=200$ and $d=100$. For the direct predictions in Table~\ref{tab:deltadetail}, each connection probability is estimated from 2,000,000 independent events under each alternative. For comparison with the simulation results, we also give a simple Gumbel approximation to the connection probabilities. Conditional on the source observation, we first approximate the squared distances by independent Gaussian variables with matching means and variances. We then use the extreme-value approximation for Gaussian minima, which has a Gumbel form. Using a common scale for the two minima gives a logistic formula for the connection probability (see the Supplement). Averaging these conditional probabilities over source observations gives the approximate connection probabilities.

Given estimates $\widehat p_1,\widehat p_2$ of the connection probabilities, we predict the component means by
\[
 \widehat\Delta_w=\frac n2(\widehat p_1+\widehat p_2-2p_{0,N})\,\widehat b_w,
 \qquad
 \widehat\Delta_d=n(\widehat p_1-\widehat p_2)\,\widehat b_d,
\]
where $\widehat b_j=500^{-1}\sum_{b=1}^{500}\sigma_j(G_b)^{-1}$ for $j\in\{w,d\}$. Here $G_b$ are 500 independent graphs generated under the same alternative, and $\sigma_j(G_b)$ is the exact permutation standard deviation of $R_j$ on graph $G_b$.

The independent simulation row averages the components over 1,000 further graphs; both graph samples are independent of the connection calculations. Under location, the simulated mean vector is $(2.649,0.032)$, with the mean shift concentrated in $Z_w$. Under scale, it is $(0.234,3.859)$, with a much larger mean shift in $Z_d$. The Gumbel approximation also captures this distinction, although its predictions for the location mean of $Z_w$ and the scale mean of $Z_d$ are below the corresponding simulation estimates.

\begin{table}[!t]
\centering
\caption{Predicted and simulated mean vectors for $(Z_w,Z_d)$ at $N=200$ and $d=100$.}
\label{tab:deltadetail}
\begin{tabular}{@{}lcc@{}}
\toprule
calculation & location $\|\mu\|=2$ & scale $\sigma=1.08$\\
\midrule
Independent simulation & $(2.649,0.032)$ & $(0.234,3.859)$\\
direct connection & ${(2.696,0.005)}$ & ${(0.237,3.865)}$\\
conditional-law MC & $(2.693,-0.040)$ & $(0.227,3.865)$\\
Gumbel connection & $(2.453,0.016)$ & $(0.198,3.417)$\\
\bottomrule
\end{tabular}
\end{table}

\subsection{Finite-sample power comparisons}
\label{sec:finitepower}

In Table~\ref{tab:power}, we compare rejection probabilities for $S=Z_w^2+Z_d^2$ at level .05, using the same $N=200$, $d=100$ models. The 500 graphs used above also estimate the sampling covariance $\widehat\Omega_{\rm MC}$ of $(Z_w,Z_d)$ under each model. Combining each predicted mean vector in Table~\ref{tab:deltadetail} with the same covariance estimate for that model gives a bivariate normal approximation. We estimate the probability that a vector from this distribution has squared norm exceeding $\chi^2_{2,0.95}$.

\begin{table}[!t]
\centering
\caption{Rejection probabilities for $S=Z_w^2+Z_d^2$ at level .05 for $N=200$ and $d=100$. }
\label{tab:power}
\begin{tabular}{@{}lccc@{}}
\toprule
calculation & null & location $\|\mu\|=2$ & scale $\sigma=1.08$\\
\midrule
Simulation: $\chi^2$ cutoff & 0.046 & 0.653 & 0.980\\
Simulation: permutation cutoff & 0.049 & 0.656 & 0.967\\
direct $p$, $\widehat\Omega_{\rm MC}$ & 0.045 & 0.655 & 0.983\\
conditional-law MC $p$, $\widehat\Omega_{\rm MC}$ & 0.047 & 0.654 & 0.982\\
Gumbel $p$, $\widehat\Omega_{\rm MC}$ & 0.046 & 0.569 & 0.932\\
\bottomrule
\end{tabular}
\end{table}

The simulation with the $\chi^2$ cutoff uses the 1,000 independent graphs described above. A separate set of 1,000 graphs per model gives the simulation results with the permutation cutoff, using 199 permutations per graph and randomized tie ranks. Thus the connection estimates, covariance estimates, and simulated rejection probabilities come from independent simulations. The null predictions use the exact connection probability $(n-1)/(2n-1)$.

The predictions based on connection probabilities estimated by direct simulation and conditional-law Monte Carlo agree closely with the independent $\chi^2$-cutoff simulation results. The Gumbel approximation gives lower power predictions under both alternatives. Under scale, the estimated covariance $\widehat\Omega_{\rm MC}$ has diagonal $(0.612,0.527)$. Holding the Gumbel mean prediction fixed and replacing this covariance by $I_2$ changes the predicted rejection probability from 0.932 to 0.876. The permutation null rejection estimate is 0.049, with Monte Carlo standard error 0.007. The Supplement gives the full computational details, Gaussian-minimum and Gumbel connection approximations, results at $d=500$, and comparisons of the null rejection probabilities at three sample sizes.

\section{Conclusion}
\label{sec:discussion}

Dimension affects scale detection through both the relative mean shifts of the two within-sample counts and their fluctuations. For the directed 1-NN test with balanced isotropic Gaussian samples, the local-power results show that scale changes are detected by the weighted component in one dimension and by the difference component in logarithmic dimension. Location changes are detected by the weighted component in both regimes.

\subsection*{Acknowledgment}
This research was supported in part by NSF grant DMS-2311399.

\bibliographystyle{imsart-nameyear}
\bibliography{references}

\end{document}